\RequirePackage{fix-cm} 
\documentclass[a4paper, twoside, reqno, dvips, 12pt]{amsart}
\usepackage{fixltx2e}   

\usepackage[latin1]{inputenc}
\usepackage[T1]{fontenc}

\usepackage{eucal}
\usepackage{esint}
\usepackage{dsfont}
\usepackage{xspace}
\usepackage{amsgen}
\usepackage{amsthm}
\usepackage{amssymb}
\usepackage{amsmath}
\usepackage{upgreek}
\usepackage{amsfonts}
\usepackage{MnSymbol}
\usepackage{mathrsfs}
\usepackage{mathtools}
\usepackage[nice]{nicefrac}

\usepackage{a4wide}

\usepackage[dvipsnames, table]{xcolor}
\definecolor{bckg}{RGB}{20.8, 20.8, 20.8}
\definecolor{oneblue}{rgb}{0.0, 0.0, 0.85}
\definecolor{Lightblue}{RGB}{214, 214, 214}
\definecolor{bluepigment}{rgb}{0.2, 0.2, 0.6}
\definecolor{charcoal}{rgb}{0.21, 0.27, 0.31}
\definecolor{denimblue}{rgb}{0.08, 0.38, 0.74}
\definecolor{Lightgray}{rgb}{0.89, 0.89, 0.89}
\definecolor{darkgrey}{rgb}{0.273, 0.281, 0.30}
\definecolor{darkelectricblue}{rgb}{0.33, 0.41, 0.47}

\usepackage[sort&compress, comma, square, numbers]{natbib}

\usepackage{psfrag}
\usepackage{graphicx}
\usepackage{subfigure}
\usepackage{morefloats}
\usepackage{indentfirst}

\usepackage[perpage, symbol]{footmisc}

\usepackage{acronym}
\usepackage{microtype}
\usepackage[labelsep=period,%
            labelfont={bf,sf,color=bluepigment},%
            justification=raggedright]{caption}

\usepackage[usenames, dvipsnames, pdf]{pstricks}
\usepackage{epsfig}
\usepackage{pst-grad} 
\usepackage{pst-plot} 

\usepackage[colorlinks,
           urlcolor=oneblue,
           linkcolor=denimblue,
           citecolor=NavyBlue,
           bookmarksopen=false,
           pdfpagemode=UseNone,
           pagebackref]{hyperref}
           
\usepackage[explicit]{titlesec}

\titleformat{\section}
  {\color{NavyBlue}\Large\sffamily\bfseries}
  {}
  {0em}
  {\colorbox{bckg!5}{\parbox{\dimexpr\linewidth-2\fboxsep\relax}{\centering\thesection. #1}}}
  [\vspace*{0.33em}]

\titleformat{name=\section,numberless}
  {\color{NavyBlue}\Large\sffamily\bfseries}
  {}
  {0.0em}
  {\colorbox{bckg!10}{\parbox{\dimexpr\linewidth-2\fboxsep\relax}{\centering#1}}}
  [\vspace*{0.33em}]

\titleformat{\subsection}
  {\color{NavyBlue}\large\sffamily\bfseries}
  {}
  {0.0em}
  {\colorbox{bckg!5}{\parbox{\dimexpr\linewidth-2\fboxsep\relax}{\centering\thesubsection. #1}}}
  [\vspace*{0.33em}]

\titleformat{name=\subsection,numberless}
  {\color{NavyBlue}\Large\sffamily\bfseries}
  {}
  {0em}
  {\colorbox{bckg!10}{\parbox{\dimexpr\linewidth-2\fboxsep\relax}{\centering#1}}}
  [\vspace*{0.33em}]

\titleformat{\subsubsection}
  {\color{bluepigment}\sffamily\normalsize\bfseries}
  {\thesubsubsection}
  {0.5em}
  {#1}
  [\vspace*{0.33em}]

\titleformat{\paragraph}[runin]
  {\color{bluepigment}\sffamily\small\bfseries}
  {}
  {0em}
  {#1}

\titlespacing{\section}{1.0em}{1.5em plus 2pt minus 2pt}%
{1.0em plus 2pt minus 2pt}[0em]
\titlespacing{\subsection}{1.0em}{1.5em plus 2pt minus 2pt}%
{1.0em}[0em]
\titlespacing{\subsubsection}{1.0em}{1.5em plus 2pt minus 2pt}%
{1.0em plus 2pt minus 2pt}[0em]

\usepackage{titletoc}

\contentsmargin{0.5em}
\newlength{\tocsep} 
\titlecontents{section}[\tocsep]
  {\addvspace{10pt}\bfseries\sffamily}
  {\contentslabel[\thecontentslabel]{\tocsep}}
  {}
  {\ \titlerule*[0.75pc]{.}\ \thecontentspage}
  []
\titlecontents{subsection}[\tocsep]
  {\addvspace{8pt}\sffamily}
  {\contentslabel[\thecontentslabel]{\tocsep}}
  {}
  {\ \titlerule*[0.5pc]{.}\ \thecontentspage}
  []
\titlecontents*{subsubsection}[\tocsep]
  {\addvspace{2pt}\footnotesize\sffamily}
  {}
  {}
  {\ \titlerule*[0.35pc]{.}\ \thecontentspage}
  [\\*]

\makeatletter
\def\@setauthors{%
  \begingroup
  \def\thanks{\protect\thanks@warning}%
  \trivlist
  \centering\footnotesize \@topsep30\p@\relax
  \advance\@topsep by -\baselineskip
  \item\relax
  \author@andify\authors
  \def\\{\protect\linebreak}%
  \textsc{\normalsize\textcolor{darkelectricblue}{\authors}}%
  \ifx\@empty\contribs
  \else
    ,\penalty-3 \space \@setcontribs
    \@closetoccontribs
  \fi
  \endtrivlist
  \endgroup
}
\def\@settitle{\begin{center}%
  \baselineskip14\p@\relax
    \bfseries
    \textsc{\Large\textcolor{charcoal}{\@title}}
  \end{center}%
}
\makeatother

\usepackage{enumitem}
\setlist[description]{%
  topsep=30pt,               
  itemsep=5pt,               
  font={\bfseries\sffamily\color{NavyBlue}}, 
}

\usepackage{fancyhdr}
\usepackage{lastpage}

\newcommand*\Title{\textcolor{bluepigment}{Supraconvergence phenomenon on non-uniform grids}}
\newcommand*\Authors{\textcolor{bluepigment}{G.~Khakimzyanov \& D.~Dutykh}}
\newcommand*{\plogo}{\textcolor{gray}{{\texttt{arXiv.org} / \textsc{hal}}}} 

\numberwithin{equation}{section}

\newtheorem{lemma}{Lemma}
\newtheorem{remark}{Remark}
\newtheorem{deff}{Definition}
\newtheorem{theorem}{Theorem}

\newcommand{\up}[1]{$^{\mathrm{\small\textsf{#1}}}$} 

\newcommand{\N}{\mathds{N}}
\newcommand{\R}{\mathds{R}}

\newcommand{\I}{\mathcal{I}}

\newcommand{\ud}{\mathrm{d}}

\newcommand{\ue}{\mathrm{e}}

\newcommand{\Ll}{\mathscr{L}}

\renewcommand{\leq}{\leqslant}
\renewcommand{\geq}{\geqslant}
\newcommand{\eps}{\varepsilon}
\renewcommand{\O}{\mathcal{O}}

\newcommand{\const}{\mathrm{const}}
\newcommand{\No}{$\mathrm{N}^\circ$}

\renewcommand{\gamma}{\boldsymbol{\upgamma}}

\newcommand*{\dprime}{{\prime\prime}\mkern-1.2mu}
\newcommand*{\trprime}{{\prime\prime\prime}\mkern-1.2mu}

\newcommand{\ie}{\emph{i.e.}~}
\newcommand{\eg}{\emph{e.g.}~}

\newcommand{\etal}{\emph{et al.}~}

\newcommand{\abs}[1]{\left|#1\right|}

\newcommand{\norm}[1]{\lVert\, #1\, \rVert}

\newcommand{\od}[2]{\frac{\mathrm{d} #1}{\mathrm{d}\/#2}}

\newcommand{\eqdef}{\mathop{\stackrel{\,\mathrm{def}}{:=}\,}}
\newcommand{\defeq}{\mathop{\stackrel{\,\mathrm{def}}{=:}\,}}

\newcommand{\half}{{\textstyle{1\over2}}}

\usepackage{acronym}
\acrodef{bvp}[BVP]{Boundary Value Problem}
\acrodef{NSWE}{Nonlinear Shallow Water Equations}

\begin{document}

\title[\Title]{On supraconvergence phenomenon for second order centered finite differences on non-uniform grids}

\author[G.~Khakimzyanov]{Gayaz Khakimzyanov}
\address{\textbf{G.~Khakimzyanov:} Institute of Computational Technologies, Siberian Branch of the Russian Academy of Sciences, Novosibirsk 630090, Russia}
\email{Khak@ict.nsc.ru}
\urladdr{http://www.ict.nsc.ru/ru/structure/Persons/ict-KhakimzyanovGS}

\author[D.~Dutykh]{Denys Dutykh$^*$}
\address{\textbf{D.~Dutykh:} LAMA, UMR 5127 CNRS, Universit\'e Savoie Mont Blanc, Campus Scientifique, 73376 Le Bourget-du-Lac Cedex, France}
\email{Denys.Dutykh@univ-savoie.fr}
\urladdr{http://www.denys-dutykh.com/}
\thanks{$^*$ Corresponding author}


\begin{titlepage}
\thispagestyle{empty} 
\noindent
{\Large Gayaz \textsc{Khakimzyanov}}\\
{\it\textcolor{gray}{Institute of Computational Technologies, Novosibirsk, Russia}}
\\[0.02\textheight]
{\Large Denys \textsc{Dutykh}}\\
{\it\textcolor{gray}{CNRS, Universit\'e Savoie Mont Blanc, France}}
\\[0.16\textheight]

\colorbox{Lightblue}{
  \parbox[t]{1.0\textwidth}{
    \centering\huge\sc
    \vspace*{1.0cm}
    
    \textcolor{bluepigment}{On supraconvergence phenomenon for second order centered finite differences on non-uniform grids}
    
    \vspace*{1.0cm}
  }
}

\vfill 

\raggedleft     
{\large \plogo} 
\end{titlepage}


\newpage
\thispagestyle{empty} 
\par\vspace*{\fill}   
\begin{flushright} 
{\textcolor{denimblue}{\textsc{Last modified:}} \today}
\end{flushright}


\newpage
\maketitle
\thispagestyle{empty}


\begin{abstract}
In the present study we consider an example of a boundary value problem for a simple second order ordinary differential equation, which may exhibit a boundary layer phenomenon depending on the value of a free parameter. To this equation we apply an adaptive numerical method on redistributed grids. We show that usual central finite differences, which are second order accurate on a uniform grid, can be substantially upgraded to the fourth order by a suitable choice of the underlying non-uniform grid. Moreover, we show also that some other choices of the nodes distributions lead to substantial degradation of the accuracy. This example is quite pedagogical and we use it only for illustrative purposes. It may serve as a guidance for more complex problems.

\bigskip
\noindent \textbf{\keywordsname:} finite differences; non-uniform grids; boundary layer; boundary value problems; supraconvergence \\

\smallskip
\noindent \textbf{MSC:} \subjclass[2010]{65N06 (primary), 76M20, 65L10 (secondary)}

\end{abstract}

\newpage
\tableofcontents
\thispagestyle{empty}


\newpage
\section{Introduction}

Boundary layer phenomena are present in many applications, in particular in Fluid Mechanics and Aerodynamics \cite{Schlichting2000}. For instance, the very successful design of the \textsc{Airbus} A320's wing is mainly due to the potential flow theory with appropriate boundary layer corrections \cite{Roe2005}. Nowadays this problem is addressed mainly with numerical techniques and it represents serious challenges.

Some numerical approaches to address the boundary layer problem have been proposed since the early 60's. Historically, probably homogeneous schemes on uniform \cite{Tikhonov1961} and non-uniform \cite{Tikhonov1962} meshes were proposed first by \textsc{Tikhonov} and \textsc{Samarskii}. Later, \textsc{Il'in} introduced the so-called exponential-fitted schemes \cite{Ilin1969, Roos1994}, which were generalized recently to finite volumes as well (see \eg \cite{Eymard2006}). The idea of \textsc{Il'in} consisted in introducing a \emph{fitting factor} into the scheme and requiring that a particular exact solution satisfies the difference equation exactly. Thanks to pioneering works of \textsc{Numerov} \cite{Numerov1927, Numerov1924}, we know that on \emph{uniform grids} it is possible to construct fourth order schemes for second order differential equations on three point stencils. However, the application of this scheme to singularly perturbed problems requires the introduction of the so-called \emph{fitting factor}, which may lead to the substantial degradation of the order of convergence. For instance, the uniform (in small parameter) first order convergence was reported in \cite[Table~2]{Phaneendra2015}. We can also mention two pioneering references where the moving grids were first applied to unsteady problems in shallow water flows \cite{Sudobicher1968} and in gas dynamics \cite{Alalykin1970}. The uniform convergence of monotone finite difference operators for singularly perturbed semi-linear equations was shown in \cite{Farrell1996}. We refer to \cite{Keller1978} for a general review of numerical methods in the boundary layer theory.

In \cite[pp. 585--586]{Roe2005} one can read:
\begin{quote}
\textit{I am convinced that it should be possible to develop a general theory of the relation between the grid, the governing equations and the specific solution being computed, but only very hazy ideas how to bring such a theory about.}
\end{quote}
Our study is a little attempt towards this research direction. An earlier attempt was undertaken in \cite{Ferreira1993}. Namely, in the present manuscript we consider a singularly perturbed linear second order elliptic ODE as a model equation which exhibits the boundary layer phenomenon. In accordance with the I.~M.~\textsc{Gelfand} principle, we took the simplest non-trivial example to illustrate our point. The goal is to propose a numerical method for such problems, which is able to solve approximately this problem with an accuracy independent of the value of the perturbation parameter \cite{Roos1994}. In the beginning we explain why the classical central finite differences on a \emph{uniform mesh} is not working in practice, even if this method is fully justified from the theoretical point of view with well-known stability and convergence properties \cite{Samarskii2001}. Then, we propose a non-uniform equidistributed grid and we show that the same central difference scheme converges with the fourth order rate on this family of successively refined grids. So, just by changing the distribution of nodes in a smart way one can gain two extra orders of the accuracy! The logarithmically-distributed grids were proposed by \textsc{Bakhvalov} \cite{Bakhvalov1969}. However, they were shown to converge inevitably with the same second order rate (see \cite{Samarskii2001} for the proof). Later supraconvergence phenomena have been studied theoretically for some elliptic boundary value problems on non-uniform grids \cite{Ferreira2006, Emmrich2006, Huang2011a}.

Non-uniform grids can be used also to compute numerically blow-up solutions as in the nonlinear Schr\"odinger equation \cite{Budd2001}. See also \cite{Budd2001a} for a general review of these techniques. There is a related idea of constructing non-uniform grids in order to preserve some or all symmetries of the continuous equation at the discrete level as well \cite{Chhay2011}. We can only regret that the authors of \cite{Chhay2011} did not study theoretically the stability and convergence of the scheme depending on symmetry preservation abilities. The same idea holds for invariants \cite{Dorodnitsyn1994, Budd2001} and asymptotics \cite{Zeldovich2002, Budd2001b}. In the present study we focus essentially on the scheme approximation order depending on the underlying (non-uniform) grid.

The phenomenon of supraconvergence of central finite difference schemes is well known and it was studied rigorously in one spatial dimension in \cite{Barbeiro2005} and the 2D case was considered in \cite{Ferreira1998}. The \emph{super-supraconvergence} reported in this manuscript is achieved by using monitoring functions, which depend on lowest order derivatives comparing to examples reported in the literature so far \cite{Degtyarev1987a, Ferreira1993}. This property greatly simplifies the implementation of grid redistribution methods. As \textsc{Strang} \& \textsc{Iserles} \cite{Strang1983} discovered the link between the stability and the stencil of a numerical scheme, here we try to understand deeper a link between scheme's convergence order and the underlying grid. In particular, we show that some thouroughly chosen nodes distributions lead to the substantial improvement of the numerical solution accuracy (for a fixed scheme). We show also that some other grid choices (appearing admissible from the first sight) may totally degrade the solution accuracy. These illustrations should serve as an indication for more complex problems.

The present manuscript is organized as follows. The \acs{bvp} under consideration is described in Section~\ref{sec:bvp}. The classical discretization is described in Section~\ref{sec:uni}, while the scheme on a general non-uniform grid is provided in Section~\ref{sec:equi}. A practical equidistribution method to construct the grids is explained in Section~\ref{sec:mesh}. A series of numerical experiments on various non-uniform grids is presented in Section~\ref{sec:num} and some theoretical insight into these results is given in Sections~\ref{sec:expla} and \ref{sec:stab}. Finally, the article is completed by outlining the main conclusions and perspectives of the present study in Section~\ref{sec:concl}.


\section{The boundary value problem}
\label{sec:bvp}

Consider the following linear \acf{bvp} for an ordinary differential equation $\Ll u\ =\ 0$ of the second degree with Dirichlet-type boundary conditions on the segment $\I\ =\ [\,0,\, \ell\,]\,$:
\begin{equation}\label{eq:bvp}
  \Ll\, u\ \eqdef\ -\,\od{^2 u}{x^2}\ +\ \lambda^2\,u\ =\ 0\,, \qquad
  u\,(0)\ =\ \ue^{-\lambda\ell}\,, \quad u\,(\ell)\ =\ 1\,,
\end{equation}
where $\lambda\ \in\ \R$ is a model parameter, which can take large (but finite) values.

\begin{figure}
  \centering
  \includegraphics[width=0.65\textwidth]{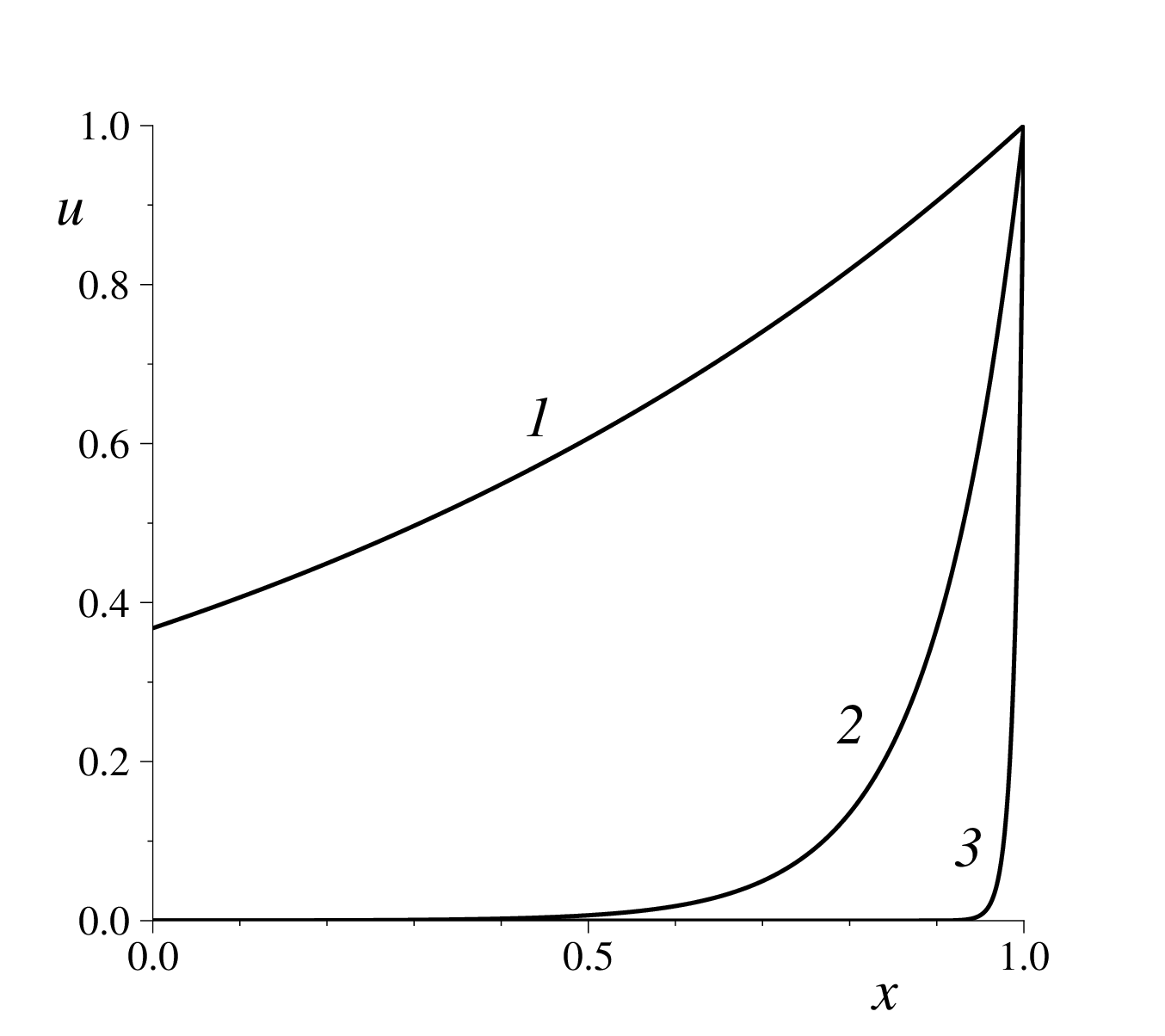}
  \caption{\small\em Exact solutions \eqref{eq:sol} for $\ell = 1$ and various values of parameter $\lambda$: (1) --- $\lambda = 1$, (2) --- $\lambda = 10$, (3) --- $\lambda = 100$.}
  \label{fig:sol}
\end{figure}

It can be readily checked that the following function of $x$ solves exactly the \acs{bvp} \eqref{eq:bvp}:
\begin{equation}\label{eq:sol}
  u\,(x)\ =\ \ue^{\,\lambda(x\ -\ \ell)}\,.
\end{equation}
However, we shall proceed as if the analytical solution \eqref{eq:sol} were not known. It will serve us only to assess the quality of a numerical solution. The peculiarity here is that for sufficiently large values of parameter $\lambda\ \gg\ 1$ the solution \eqref{eq:sol} shows a boundary layer type behaviour in the vicinity of the point $x\ =\ \ell\,$. It is illustrated in Figure~\ref{fig:sol}. Similar phenomena occur in Fluid Mechanics where they are of capital importance \eg in Aerodynamics \cite{Schlichting2000}. It justifies the choice of the problem \eqref{eq:bvp} in our study.

Usually the problem \eqref{eq:bvp} is rewritten in the literature by introducing a small parameter $\eps\ \eqdef\ 1 / \lambda^2\,$:
\begin{equation*}
  -\eps\,\od{^2 u}{x^2}\ +\ u\ =\ 0\,.
\end{equation*}
Thus, we have a singularly-perturbed \textsc{Sturm}--\textsc{Liouville} problem \cite{Nikitin1999}. 

\subsection{Discretization on a uniform grid}
\label{sec:uni}

Consider a uniform discretization of the segment $\I$ in $N$ equal segments $\I_{\,h}$ with boundaries located at $\{x_j\, =\, j\,\Delta x\}_{j\,=\,0}^{\,N}$, $\Delta x\, =\, \frac{\ell}{N}\,$. The finite difference analogue $\Ll^{\,h}\, u_{\,h}\ =\ 0$ of differential equation \eqref{eq:bvp} is
\begin{equation}\label{eq:uni}
  \Ll^{\,h}_{\,j}\{u_{\,j}\}\ =\ -\frac{u_{\,j+1}\ -\ 2\,u_{\,j}\ +\ u_{\,j-1}}{\Delta x^2}\ +\ \lambda^2\;u_{\,j}\ =\ 0\,, \qquad j\ =\ 1,\,\ldots,\,N\,-\,1\,,
\end{equation}
together with \textsc{Dirichlet}-type boundary conditions:
\begin{equation*}
  u_{\,0}\ =\ \ue^{-\lambda\ell}\,, \qquad u_{\,N}\ =\ 1\,.
\end{equation*}
It is well known that this scheme has the second order accuracy, as it follows from the local consistency error analysis:
\begin{equation*}
  e_{\,j}^{\,h}\ \eqdef\ \Ll^{\,h}_{\,j}\{u\,(x_{\,j})\}\ =\ u_{\,xxxx}\,(x_{\,j})\,\frac{\Delta x^2}{12}\ +\ \O(\Delta x^4)\ =\ \lambda^4\ue^{\lambda(x_{\,j}\ -\ \ell)}\,\frac{\Delta x^2}{12}\ +\ \O(\Delta x^4)\,.
\end{equation*}
From the last formula we can already draw some preliminary conclusions:
\begin{itemize}
  \item The proportionality constant grows as the fourth power of the parameter $\lambda$, which can take potentially large values in practically important situations.
  \item The ratio between the consistency error in the vicinity of $x\ =\ \ell$ and $x\ =\ 0$ is $\approx\exp(\lambda\ell)\,$.
\end{itemize}

We would like to mention that the stability proof of scheme \eqref{eq:uni} can be found \eg in \cite{Godunov1987}. So, according to the \textsc{Lax}--\textsc{Richtmyer} equivalence theorem \cite{Lax1956}, the scheme \eqref{eq:uni} is convergent as $\Delta x\ \to\ 0$ and the convergence rate is equal to the approximation order (\ie two in this particular case). However, despite all these good properties the scheme \eqref{eq:uni} is not usable in practice because of two practical drawbacks mentioned above (they are all related to the asymptotic limit $\lambda\to\infty$ and the boundary layer phenomenon). Please, note however that the convergence is established for a fixed value of the parameter $\lambda$.

\subsection{Non-uniform adaptive grids}
\label{sec:equi}

In order to cope with the shortcomings mentioned above, we turn to non-uniform grids by preserving the simplicity of the second-order central discretization \eqref{eq:uni}. Let $Q\ \eqdef\ [\,0,\,1\,]$ be the reference domain and consider a bijective mapping from $Q$ to $\I\,$:
\begin{equation}\label{eq:map}
  x\,(q)\,:\ q\in Q\ \mapsto\ \I\,.
\end{equation}
We require that the boundary points map into each other:
\begin{equation*}
  x\,(0)\ =\ 0\,, \qquad x\,(1)\ =\ \ell\,.
\end{equation*}
We shall also assume that the Jacobian $J(q)$ of mapping \eqref{eq:map} is bounded from below and above by some positive constants:
\begin{equation}\label{eq:jac}
  0\ <\ J_m\ \leq\ J(q)\ \eqdef\ \od{x(q)}{q}\ \leq\ J_M\ <\ \infty\,, \qquad \forall q\ \in\ Q\,.
\end{equation}
The reference domain can be discretized into $N$ equal elements $Q_{\,h}$ by nodes $\{q_j\, =\, j h\}_{\,j\,=\,0}^{\,N}\,$, $h\ =\ 1/N\,$. Strictly speaking, for our numerical purposes we can be satisfied with the knowledge of the discrete mapping $x^h\,:\ Q_{\,h}\ \mapsto\ \I_{\,h}\,$. From condition \eqref{eq:jac} (more precisely from this part: $0\ <\ J_{\,m}\ \leq\ J\,(q)$) follows that the steps of the non-uniform mesh $\I_{\,h}$ are necessarily positive, \ie
\begin{equation}\label{eq:hj}
  h_{\,j + 1/2}\ \eqdef\ x_{\,j+1}\ -\ x_{\,j}\ >\ 0\,, \qquad
  j\ =\ 0,\,1,\,\ldots,\,N\,-\,1\,.
\end{equation}
From the condition $J\,(q)\ \leq\ J_{\,M}\ <\ \infty$ follows that
\begin{equation*}
  h_{\,\max}\ \eqdef\ \max_{j\,=\,0,\,\ldots,\,N\,-\,1} h_{\,j+1/2}\ \leq\ J_{\,M}\,h\ \to\ 0\,, \qquad h\ \to\ 0\,.
\end{equation*}
We note also that if the \textsc{Jacobian} $J\,(q)\ \equiv\ \const$ (or equivalently the mapping $x\,(q)$ is linear), the resulting mesh $\I_h$ is uniform. This case is not very interesting from the adaptivity point of view. Any other map $x\,(q)\,$, which satisfies the boundary conditions with \eqref{eq:jac} defines a valid non-uniform grid $\I_{\,h}\,$. However, not all of them are equally interesting from the numerical point of view. Some practical methods to compute solution-adapted meshes will be discussed below in Section~\ref{sec:mesh}.

Equation \eqref{eq:bvp} can be posed on the reference domain $Q\,$:
\begin{equation}\label{eq:trans}
  -\frac{1}{J}\;\od{}{q}\;\biggl(\frac{1}{J}\;\od{v}{q}\biggr)\ +\ \lambda^2\, v\ =\ 0\,,
\end{equation}
with the same boundary conditions $v\,(0)\ =\ \ue^{-\,\lambda\,\ell}\,$, $v\,(1)\ =\ 1\,$. Here $v$ is the composed function $v\,(q)\ \eqdef\ u\circ x\ =\ u\,\bigl(x\,(q)\bigr)\,$. Now equation \eqref{eq:trans} is posed on domain $Q$ and, thus, it can be discretized on the \emph{uniform} grid $Q_{\,h}$ as follows
\begin{equation}\label{eq:discr}
  -\,\frac{1}{h\,J_{\,j}}\;\biggl[\frac{v_{\,j+1}\ -\ v_{\,j}}{h\,J_{\,j+1/2}}\ -\ \frac{v_{\,j}\ -\ v_{\,j-1}}{h J_{\,j-1/2}}\biggr]\ +\ \lambda^2\,v_{\,j}\ =\ 0\,, 
  \qquad j\ =\ 1,\,2,\,\ldots,\,N\,-\,1\,,
\end{equation}
where the \textsc{Jacobian} $J$ is computed as
\begin{equation*}
  J_{\,j+1/2}\ \eqdef\ \frac{x_{\,j+1}\ -\ x_{\,j}}{h}\,, \quad
  J_{\,j-1/2}\ \eqdef\ \frac{x_{\,j}\ -\ x_{\,j-1}}{h}\,, \quad
  J_{\,j}\ \eqdef\ \frac{J_{\,j-1/2}\ +\ J_{\,j+1/2}}{2}\,.
\end{equation*}
Difference equations \eqref{eq:discr} have to be completed with corresponding boundary conditions
\begin{equation*}
  v_{\,0}\ =\ \ue^{-\,\lambda\,\ell}\,, \qquad v_{\,N}\ =\ 1\,.
\end{equation*}
The last scheme \eqref{eq:discr} can be equivalently rewritten on the non-uniform grid $\I_{\,h}\,$:
\begin{equation}\label{eq:nonuni}
  -\frac{1}{h_{\,j}}\;\biggl[\,\frac{u_{\,j+1}\ -\ u_{\,j}}{h_{\,j+1/2}}\ -\ \frac{u_{\,j}\ -\ u_{\,j-1}}{h_{\,j-1/2}}\,\biggr]\ +\ \lambda^2\,u_{\,j}\ =\ 0\,, 
  \qquad j\ =\ 1,\,2,\, \ldots,\, N\, -\, 1\,,
\end{equation}
where $h_{\,j\pm 1/2}$ were defined in \eqref{eq:hj} and 
\begin{equation*}
  h_{\,j}\ \eqdef\ \frac{h_{\,j-1/2}\ +\ h_{\,j+1/2}}{2}\,.
\end{equation*}


\subsection{Adaptive mesh generation}\label{sec:mesh}

The finite difference scheme on the transformed uniform $Q_{\,h}$ and non-uniform $\I_{\,h}$ grids were formulated in the previous Section. We choose the equidistribution method presented below. A similar nodes redistribution method was explained in details for time-dependent hyperbolic problems in a companion paper \cite{Khakimzyanov2015a}. The exposition below is simpler since we deal with a stationary problem. Thus, we do not have here an additional complication of nodes motion in space and in time, which was addressed in \cite{Khakimzyanov2015a}.

A non-uniform grid $\I_h$ is given if we construct somehow the mapping $x\,(q):\; Q \mapsto \I$ and evaluate it in the nodes of the uniform grid, \ie $\I_{\,h}\ =\ x\,(Q_{\,h})\,$. In the equidistribution method it is proposed that the desired mapping $x\,(q)$ is obtained as a solution to the following nonlinear elliptic problem
\begin{equation}\label{eq:elliptic}
  \od{}{q}\;\biggl[\omega\,(x)\,\od{x}{q}\biggr]\ =\ 0\,, \qquad x\,(0)\ =\ 0\,, \quad x\,(1)\ =\ \ell\,,
\end{equation}
where $\omega\,(x)$ is the so-called \emph{monitor} function. In order to have a well-posed problem \eqref{eq:elliptic}, the function $\omega\,(x)$ has to be sufficiently smooth and positive valued. The equidistribution principle can be readily obtained by integrating \eqref{eq:elliptic} once in $q-$space
\begin{equation}\label{eq:equi_cont}
  \omega\,\bigl(x(q)\bigr)\, J(q)\ \equiv\ C\ =\ \const, \qquad \forall q\ \in\ Q,
\end{equation}
and another time in $x$-space on the element $[x_{\,j},\, x_{\,j+1}]\,$:
\begin{equation}\label{eq:cond}
  \int_{x_{\,j}}^{x_{\,j+1}} \omega\,(x)\;\ud x\ =\ C\,h\ =\ \const\,.
\end{equation}
The last identity explains the name of the equidistribution principle, \ie the quantity $\omega\,(x)$ is distributed uniformly over the cells $\I_{\,h}\,$. In other words, in the areas where $\omega\,(x)$ takes high values, the elements $[\,x_{\,j},\, x_{\,j+1}]$ are smaller in order to satisfy the condition \eqref{eq:cond}. The constant $C$ is not arbitrary and it can be computed exactly:
\begin{equation*}
  C\ =\ \int_{0}^1 \omega\,\bigl(x\,(q)\bigr)\, J\,(q)\,\ud q\ =\ \int_0^\ell \omega\,(x)\,\ud x\,.
\end{equation*}

In practice, one has to choose the monitor function $\omega\,(x)$ and solve the nonlinear elliptic Boundary Value Problem (BVP) \eqref{eq:elliptic} in order to obtain the required mapping $x\,(q)\,$. Since we need to know only this mapping in the grid nodes $q_{\,j}\,$, the problem \eqref{eq:elliptic} can be discretized using the central finite differences as well
\begin{equation}\label{eq:nbvp}
  \frac{1}{h}\;\Bigl[\,\omega(x_{j+1/2})\;\frac{x_{j+1}\ -\ x_{j}}{h}\ -\ \omega(x_{j-1/2})\;\frac{x_{j}\ -\ x_{j-1}}{h}\,\Bigr]\ =\ 0\,,
  \qquad j\ =\ 1,\,\ldots,\, N\,-\,1\,,
\end{equation}
with discrete boundary conditions $x_{\,0}\ =\ 0\,$, $x_{\,N}\ =\ \ell\,$. Equations \eqref{eq:nbvp} are solved iteratively using the following linearization:
\begin{equation*}
  \frac{1}{h}\;\biggl[\,\omega(x_{\,j+1/2}^{\,(n)})\;\frac{x_{\,j+1}^{\,(n+1)}\ -\ x_{\,j}^{\,(n+1)}}{h}\ -\ \omega(x_{\,j-1/2}^{\,(n)})\;\frac{x_{\,j}^{\,(n+1)}\ -\ x_{\,j-1}^{\,(n+1)}}{h}\,\biggr]\ =\ 0\,, \qquad n\ \in\ \N\,.
\end{equation*}
The iterations are continued until the convergence within a prescribed tolerance is achieved. The linearized equations can be solved efficiently using the tridiagonal matrix algorithm \cite{Fedorenko1994}.

A solution $\{x_{\,j}\}_{\,j\,=\,0}^{\,N}$ to the nonlinear BVP \eqref{eq:nbvp} satisfies the discrete version of the equidistribution principle \eqref{eq:cond}:
\begin{equation*}
  \omega\,(x_{\,j+1/2})\;\frac{h_{\,j+1/2}}{h}\ =\ C_{\,h}\ =\ \const\,, \qquad j\ =\ 0,\,1,\,\ldots,\,N\,-\,1\,,
\end{equation*}
or equivalently
\begin{equation*}
  \omega\,(x_{\,j+1/2})\;J_{\,j+1/2}\ =\ C_{\,h}\ =\ \const\,, \qquad j\ =\ 0,\,1,\,\ldots,\, N\,-\,1\,.
\end{equation*}
The last identity is called the \emph{discrete equidistribution principle}.

\subsubsection{Example}

In order to illustrate the use of the equidistribution principle in practice, we make the following choice of the monitor function
\begin{equation}\label{eq:example}
  \omega\,(x)\ =\ \bigl(u_{xxx}\bigr)^{\,\frac{1}{4}}\ \propto\ \bigl(u_{x}\bigr)^{\frac{1}{4}}\,, \qquad x\ \in\ \I\,,
\end{equation}
where $u\,(x)$ is the boundary layer solution \eqref{eq:sol}. The proportionality $u_{\,xxx} \propto u_{\,x}$ follows from equation \eqref{eq:bvp} and we can replace the monitor function by a proportional one without changing the grid due to the following
\begin{lemma}\label{lem:prop}
Coordinate transformations $x\ =\ x_{\,1}\,(q)$ and $x\ =\ x_{\,2}\,(q)$ which are solutions to problem \eqref{eq:elliptic} with two monitoring functions $\omega\ =\ \omega_{\,1}\,(x)$ and $\omega\ =\ \omega_{\,2}\,(x)\ \equiv C\,\omega_{\,1}\,(x)\,$, where $C\ \neq\ 0\,$, generate the same mesh.
\end{lemma}
\begin{proof}
Trivial.
\end{proof}

In the case of the monitoring function given in \eqref{eq:example}, the nonlinear BVP \eqref{eq:elliptic} can be solved exactly by using the equidistribution principle \eqref{eq:equi_cont}, with constant
\begin{equation*}
  C\ =\ \int_0^\ell \lambda^{\,\frac{3}{4}}\ue^{\frac{\lambda\,(x\ -\ \ell)}{4}}\;\ud x\ =\ \frac{4}{\lambda^{\,\frac{1}{4}}}\bigl(1\ -\ \ue^{-\,\frac{\lambda\ell}{4}}\bigr)\,.
\end{equation*}
Consequently, \eqref{eq:equi_cont} reads
\begin{equation*}
  \lambda^{\,\frac{3}{4}}\ue^{\frac{\lambda(x\ -\ \ell)}{4}}\;\od{x}{q} \ =\ \frac{4}{\lambda^{\,\frac{1}{4}}}\bigl(1\ -\ \ue^{-\frac{\lambda\ell}{4}}\bigr)\,.
\end{equation*}
From the last equation the mapping $x\,(q)$ can be found exactly
\begin{equation*}
  x\,(q)\ =\ \ell\ +\ \frac{4}{\lambda}\;\ln\bigl(\,q\ +\ (1\ -\ q)\,\ue^{-\,\frac{\lambda\ell}{4}}\,\bigr)\,.
\end{equation*}
Below it will become clear why we pay a special attention to this particular distribution of the nodes.

\subsubsection{Generalizations}

Example shown above can be easily generalized by taking the following monitor function:
\begin{equation}\label{eq:beta}
  \omega_{\beta}\,(x)\ =\ \bigl(u_{\,x}\,(x)\bigr)^{\,\beta}\,, \qquad \beta\in\R_0^+\,, \qquad x\ \in\ \I\,.
\end{equation}
By performing the same computations as in the example above, one can show that this monitor function $\omega_{\,\beta}\,(x)$ yields the following mapping between the reference domain $Q$ and computational domain $\I\,$:
\begin{equation}\label{eq:grids}
  x\,(q)\ =\ \left\{\begin{array}{cl}
    \ell\ +\ \frac{1}{\beta\lambda}\;\ln\bigl[\,q\ +\ (1\ -\ q)\,\ue^{-\beta\lambda\ell}\,\bigr], & \beta\ \neq\ 0\,, \\
    q\ell, & \beta\ =\ 0\,.
  \end{array}\right.
\end{equation}
For instance, one can see that the case $\beta\ =\ 0$ yields a uniform distribution of the nodes, which is not very interesting in view of generating adaptive redistributed meshes. Below we shall consider also a particular case when $\beta = \frac{1}{2}\,$, which gives the mapping $x\,(q)$ given by
\begin{equation*}
  x\,(q)\ =\ \ell\ +\ \frac{2}{\lambda}\;\ln\bigl[\,q\ +\ (1\ -\ q)\,\ue^{-\,\frac{\lambda\ell}{2}\,}\bigr]\,.
\end{equation*}


\section{Numerical results}
\label{sec:num}

In order to measure the quality of the numerical solution we compute its `distance' to the reference solution given by \eqref{eq:sol}:
\begin{equation*}
  ||\eps_h||_\infty\ \equiv\ ||u_h\ -\ (u)_h||_\infty\ =\ \max_{0 \leq j\leq N} |u_j\ -\ u(x_j)|\,.
\end{equation*}
The convergence order $p$ of a scheme can be estimated numerically as
\begin{equation*}
  p\ =\ \log_2\;\biggl[\,\frac{||\eps_{h}||_\infty}{||\eps_{h/2}||_\infty}\biggr]\,.
\end{equation*}
The parameter $p$ has to be computed on a sequence of refined meshes to have a more robust estimation.

The numerical results are presented in Table~\ref{tab:conv} for various choices of the monitor function of the general form \eqref{eq:beta}. The corresponding non-uniform meshes $\{x_j\ =\ x\,(q_{\,j})\}_{\,j\,=\,0}^{\,N}$ were computed analytically above in \eqref{eq:grids}. The choice of the monitor function $\omega\,(x)\ \equiv\ 1$ (\ie $\beta\ =\ 0$) corresponds to a uniform grid and in this case we recover the theoretical 2\up{nd} order convergence. However, the main focus of this study is on non-uniform grids. The most surprising result is the performance of the monitor function $\omega\,(x)\ =\ \bigl(u_{\,x}\,(x)\bigr)^{\,1/4}\,$. In this case we observe a fair 4\up{th} order convergence! When the monitor function is changed to $\omega\,(x)\ =\ \bigl(u_{\,x}\,(x)\bigr)^{\,1/2}$ we come back to the 2\up{nd} order convergence, even if the error norm $\|\varepsilon_h\|_{\infty}$ is approximately $10$ times lower. However, the most catastrophic results are observed for $\omega\,(x)\ =\ \bigl(u_{\,x}\,(x)\bigr)^{\,2}\,$, since the convergence order falls down to $p\ =\ \half\,$. This distribution of nodes is certainly to be avoided in practical numerical simulations.

\begin{table}
\centering
\begin{tabular}{|c|c|c|c|c|c|c|c|c|}
  \hline\hline
  $N$  &\multicolumn{2}{|c|}{$\omega\ =\ 1$}
  &\multicolumn{2}{|c|}{$\omega\ =\ \bigl(u_{\,x}\,(x)\bigr)^{\,1/4}$}
  & \multicolumn{2}{|c|}{$\omega\ =\ \bigl(u_{\,x}\,(x)\bigr)^{\,1/2}$}
  & \multicolumn{2}{|c|}{$\omega\ =\ \bigl(u_{\,x}\,(x)\bigr)^{\,2}$} \\ \cline{2-9}
  & $\|\varepsilon_h\|_{\infty}$ & $p$
  & $\|\varepsilon_h\|_{\infty}$ & $p$
  & $\|\varepsilon_h\|_{\infty}$ & $p$
  & $\|\varepsilon_h\|_{\infty}$ & $p$  \\ \hline
  10 & $0.141\cdot 10^{-1}$ & --- &
  $0.146\cdot 10^{-4}$ &---&
  $0.456\cdot 10^{-2}$ &---&
  $0.193$              &--- \\
  20 & $0.375\cdot 10^{-2}$ & 1.93 &
  $0.883\cdot 10^{-6}$&4.04&
  $0.101\cdot 10^{-2}$&2.17&
  $0.137$ & 0.49 \\
  40 & $0.953\cdot 10^{-3}$ & 1.96 &
  $0.548\cdot 10^{-7}$&4.01&
  $0.220\cdot 10^{-3}$&2.20&
  $0.960\cdot 10^{-1}$ & 0.49\\
  80 & $0.239\cdot 10^{-3}$ & 2.0 &
  $0.342\cdot 10^{-8}$ &4.0 &
  $0.512\cdot 10^{-4}$ &2.10 &
  $0.668\cdot 10^{-1}$ &0.49 \\
  160 &  $0.599\cdot 10^{-4}$ & 2.0 &
  $0.214\cdot 10^{-9}$ & 4.0&
  $0.127\cdot 10^{-4}$ & 2.0&
  $0.463\cdot 10^{-1}$ & 0.49\\
  320 &  $0.150\cdot 10^{-4}$ & 2.0 &
  $0.136\cdot 10^{-10}$ &3.97 &
  $0.317\cdot 10^{-5}$ &2.0 &
  $0.319\cdot 10^{-1}$ &0.59 \\
  640 &  $0.374\cdot 10^{-5}$ & 2.0 &
  $0.836\cdot 10^{-12}$ &4.03 &
  $0.792\cdot 10^{-6}$ &2.0 &
  $0.220\cdot 10^{-1}$ &0.59 \\
  \hline\hline
\end{tabular}
\bigskip
\caption{\small\em Numerical estimation of the convergence order of the scheme \eqref{eq:discr} for different choices of the monitor function $\omega\,(x)\,$. The computations are performed for $\ell\ =\ 1$ and $\lambda\ =\ 10\,$.}
\label{tab:conv}
\end{table}

\begin{figure}
  \centering
  \includegraphics[width=0.99\textwidth]{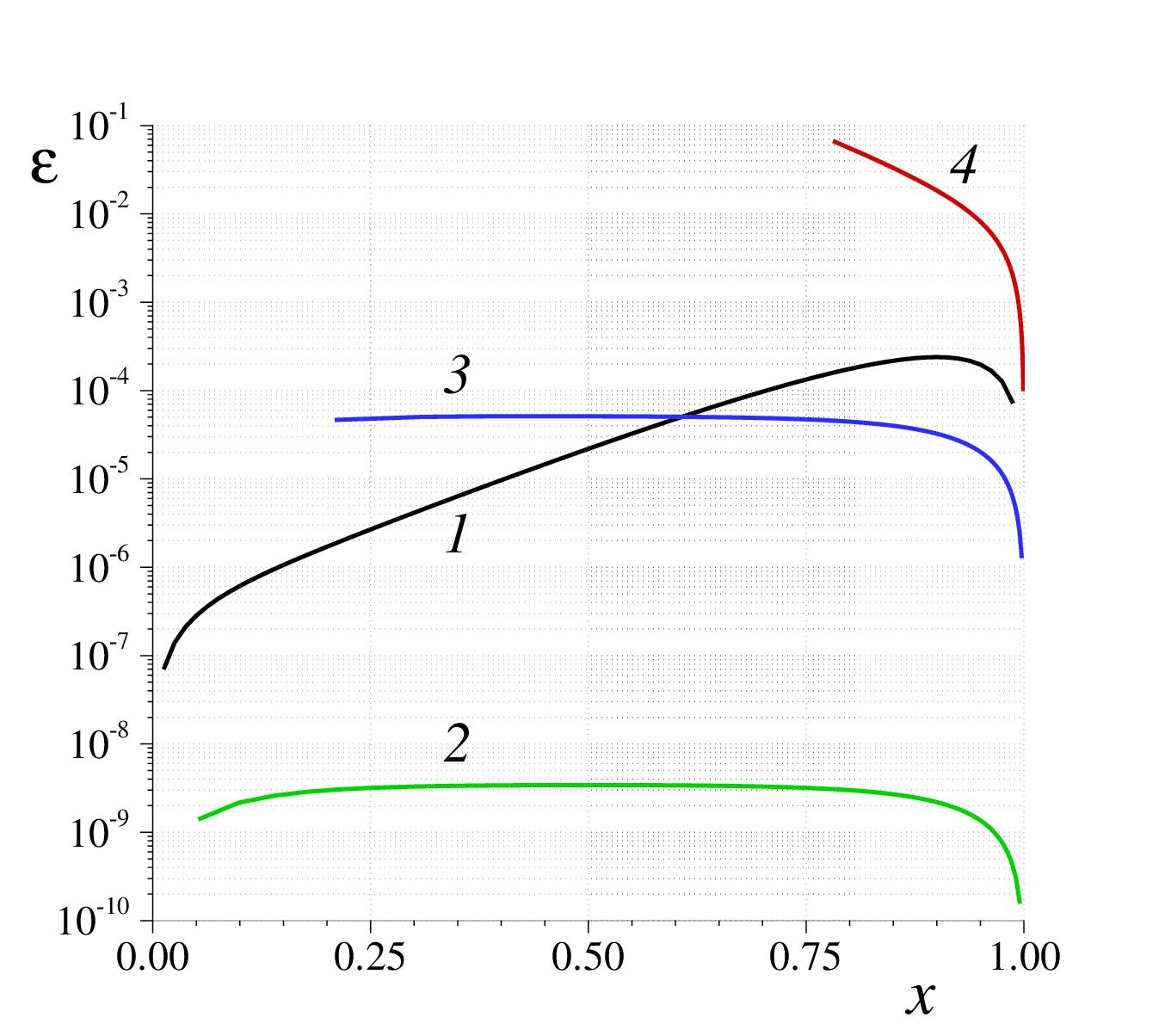}
  \caption{\small\em Distribution of the numerical error in the computational domain for various choices of the monitor function $\omega\,(x)\,$: (1) --- $\omega\,(x)\ \equiv\ 1$ (\ie the uniform grid, (2) --- $\omega\,(x)\ =\ \bigl(u_{\,x}\,(x)\bigr)^{\,\frac{1}{4}}$ (optimal case), (3) --- $\omega\,(x)\ =\ \bigl(u_{\,x}\,(x)\bigr)^{\,\frac{1}{2}}\,$, (4) --- $\omega\,(x)\ =\ \bigl(u_{\,x}\,(x)\bigr)^{\,2}\,$. The problem parameters are $\ell\ =\ 1\,$, $\lambda\ =\ 10$ and only $N\ =\ 80$ nodes are used.}
  \label{fig:const}
\end{figure}

It is important to see how the error is distributed along the computational domain $\I\,$. For the same choices of the monitor function $\omega_\beta\,(x)$ it is depicted in Figure~\ref{fig:const}. One can see that the uniform grid leads to exponentially inaccurate results in the boundary layer (\ie the vicinity of $x\ =\ \ell$), despite its good theoretical properties on the paper. Without any surprise the supraconvergent case of $\omega_{1/4}\,(x)$ leads the lowest absolute error. However, it is surprising to see that $\omega_{1/2}\,(x)$ suppresses the main drawback of the uniform mesh (while preserving its order of convergence) --- the error does not explode anymore in the boundary layer. Finally, one can see that $\omega_{2}(x)$ is a poor choice since it leads to huge errors on the left end of the computational domain (\ie $x\ =\ 0$). What happens is that we put too many nodes in the boundary layer and we forget about the rest of the domain. Without any surprise this strategy cannot lead to good numerical results.

In order to complete the numerical study, we slightly change the monitor function to make it closer to what is actually used in practice \cite{Khakimzyanov2015a}. Namely, we take the following family of monitor functions parametrized by two positive real parameters $\alpha, \beta \in \R_0^+$:
\begin{equation}\label{eq:moni}
  \omega_{\,\alpha,\,\beta}\,(x)\ =\ 1\ +\ \alpha\;\abs{\,u_{\,x}}^{\,\beta}\,,
\end{equation}
where this time $u(x)$ is \emph{not} the exact solution given in \eqref{eq:sol}, but the solution being computed numerically. Thus, the grid generation process becomes iterative and the stopping criterium is 
\begin{equation*}
  \norm{u_{\,h}^{\,n+1}\ -\ u_{\,h}^{\,n}}_{\,\infty}\ <\ \eps\,,
\end{equation*}
where $\eps$ is a tolerance parameter taken to be $\eps = 10^{-10}$ in the simulations performed below. The iterative procedure is explained above in Section~\ref{sec:mesh}. A grid generated in this way will be called an \emph{adaptive grid} to the solution with the monitor function $\omega_{\alpha,\beta}\,(x)$ defined in \eqref{eq:moni}. The numerical results reported in Table~\ref{tab:res} show that here again the lowest errors in the solution are achieved for the parameter $\beta\ =\ 1/4\,$. For $\beta\ =\ 1/2$ (as well as for $\beta\ =\ 1$ and $\beta\ =\ 2$) we can see that the error does not depend \emph{monotonically} on the other parameter $\alpha\,$. The results for $\beta\ =\ 2$ are even worse than for the corresponding uniform grid.

In general, the optimal determination of the monitoring function parameters $\alpha\,$, $\beta$ is a complex problem. In principle, one could investigate the numerical performance of even more general monitor functions such as
\begin{equation*}
  \omega\,(x)\ =\ 1\ +\ \alpha_0\;|u|^{\,\beta_0}\ +\ \alpha_1\;|u_x|^{\,\beta_1}\ +\ \alpha_2\;|u_{xx}|^{\,\beta_2}\,,
\end{equation*}
but it is beyond the scope of the present study.

\begin{table}
\centering
\begin{tabular}{|c|c|c|c|c|c|c|c|c|c|c|}
  \hline\hline 
  $\alpha$ &\multicolumn{2}{|c|}{$\beta = 1/8$}
  & \multicolumn{2}{|c|}{$\beta = 1/4$}
  & \multicolumn{2}{|c|}{$\beta = 1/2$}
  & \multicolumn{2}{|c|}{$\beta = 1$}
  & \multicolumn{2}{|c|}{$\beta = 2$}\\ \cline{2-11}
  & $\|\varepsilon_h\|_{\infty}$ & $n$
  & $\|\varepsilon_h\|_{\infty}$ & $n$
  & $\|\varepsilon_h\|_{\infty}$ & $n$
  & $\|\varepsilon_h\|_{\infty}$ & $n$ 
  & $\|\varepsilon_h\|_{\infty}$ & $n$  \\ \hline
  0 & $0.375\cdot 10^{-2}$ & 1 &
  $0.375\cdot 10^{-2}$ & 1 &
  $0.375\cdot 10^{-2}$ & 1 &
  $0.375\cdot 10^{-2}$ & 1 & 
  $0.375\cdot 10^{-2}$ & 1 \\
  0.1 & $0.330\cdot 10^{-2}$ & 6 & $0.288\cdot 10^{-2}$ & 7 &
  $0.206\cdot 10^{-2}$& 8 &
  $0.715\cdot 10^{-3}$ & 12 &
  $0.382\cdot 10^{-2}$ & 35 \\
  0.5 & $0.230\cdot 10^{-2}$ & 8 & $0.141\cdot 10^{-2}$& 10 &
  $0.358\cdot 10^{-3}$& 13 &
  $0.150\cdot 10^{-2}$ & 22 &
  $0.135\cdot 10^{-1}$ & 92 \\
  1 & $0.182\cdot 10^{-2}$ & 10 & $0.816\cdot 10^{-3}$ & 12 &
  $0.321\cdot 10^{-3}$ & 16 &
  $0.176\cdot 10^{-2}$ & 30 &
  $0.377\cdot 10^{-1}$ & 167 \\
  2 & $0.142\cdot 10^{-2}$ & 10 & $0.423\cdot 10^{-3}$ & 15 &
  $0.483\cdot 10^{-3}$ & 22 &
  $0.230\cdot 10^{-2}$ & 42 &
  $0.841\cdot 10^{-1}$ & 699 \\
  $10$ & $0.951\cdot 10^{-3}$ & 13 & $0.824\cdot 10^{-4}$ & 20 &
  $0.630\cdot 10^{-3}$ & 38 &
  $0.750\cdot 10^{-2}$ & 123 &
  $0.227$ &73 \\
  $10^2$ & $0.832\cdot 10^{-3}$ & 14 & $0.854\cdot 10^{-5}$ & 22 &
  $0.827\cdot 10^{-3}$ & 46 &
  $0.630\cdot 10^{-1}$ & 45 &
  $0.204$ &66 \\
  $10^3$ & $0,820\cdot 10^{-3}$ & 14 & $0,132\cdot 10^{-5}$ & 22 &
  $0.113\cdot 10^{-2}$ & 46 &
  $0.743\cdot 10^{-1}$ & 36 &
  $0.202$ & 65 \\
  $10^4$ & $0.819\cdot 10^{-3}$ & 14 & $0.644\cdot 10^{-6}$ & 23 &
  $0.117\cdot 10^{-2}$ & 46 &
  $0.754\cdot 10^{-1}$ & 37 &
  $0.202$ & 64 \\
  \hline\hline
\end{tabular}
\bigskip
\caption{\small\em Effect of monitoring function \eqref{eq:moni} parameters $\alpha\,$, $\beta$ on the accuracy of the fully converged numerical solution. Here we report also the number $n$ of iterations needed to achieve the prescribed tolerance $\eps = 10^{-10}\,$. Problem parameters are $\ell = 1\,$, $\lambda = 10$ and only $N = 20$ nodes are used.}
\label{tab:res}
\end{table}


\section{Explanations}
\label{sec:expla}

In order to have some theoretical insight into the numerical results presented above (especially we would like to explain the surprising performance of the ``magic'' choice $\beta\ =\ 1/4$), we study the approximation properties of the scheme \eqref{eq:nonuni} on general non-uniform grids:
\begin{equation*}
  \psi_{\,j}\ =\ -\,\frac{1}{h_{\,j}}\;\biggl[\frac{u\,(x_{\,j+1})\ -\ u\,(x_{\,j})}{h_{\,j+1/2}}\ -\ \frac{u\,(x_{\,j})\ -\ u\,(x_{\,j-1})}{h_{\,j-1/2}}\biggr]\ +\ \lambda^2\; u\,(x_{\,j})\,.
\end{equation*}
The smallness of $\psi_{\,j}$ in terms of $h$ (the uniform step in the reference domain) will characterize the consistency error of the scheme \eqref{eq:nonuni}. We assume the following
\begin{deff}
  The difference scheme $\Ll^{\,h}\, u_{\,h}\ =\ 0$ approximates the problem $\Ll\, u\ =\ 0$ on its solution with the order $p\,$, if there exists a constant $C\ =\ \const$ and a natural number $N_0\ \in\ \N$ such that $\forall N\ \geq\ N_0$
  \begin{equation*}
    \norm{\,\psi_{\,h}\,}_\infty\ \leq\ C\, h^{\,p}\,.
  \end{equation*}
\end{deff}

Let us compute the consistency error $\psi_{\,j}$ for the scheme \eqref{eq:nonuni}. To achieve this goal we develop the solution into local Taylor expansions and composing finite differences which appear in \eqref{eq:nonuni}:
\begin{align*}
  \frac{u\,(x_{\,j+1})\ -\ u\,(x_{\,j})}{h_{\,j+1/2}}\ &=\ u'\ +\ \frac{h_+}{2}\;u''\ +\ \frac{h_+^2}{6}\;-u'''\ +\ \frac{h_+^3}{24}\; u^{(4)}\ +\ \frac{h_+^4}{120}\; u^{(5)}\ +\ \O\,(h_+^5)\,, \\
  \frac{u\,(x_{\,j})\ -\ u\,(x_{\,j-1})}{h_{\,j-1/2}}\ &=\ u'\ -\ \frac{h_-}{2}\; u''\ +\ \frac{h_-^2}{6}\; u'''\ -\ \frac{h_-^3}{24}\; u^{(4)}\ +\ \frac{h_-^4}{120}\; u^{(5)}\ +\ \O(h_-^5)\,,
\end{align*}
where $h_{\pm}\ \eqdef\ h_{\,j \pm 1/2}$ and the derivatives on the right hand side are evaluated at $x\ =\ x_{\,j}\,$. We did not specify it for the sake of compactness, but this assumption will hold below. Now we can estimate asymptotically
\begin{multline}\label{eq:psij}
  \psi_{\,j}\ =\ -u''\ +\ \lambda^2\,u\ -\ \frac{h_+ - h_-}{3}\;u'''\ -\ \frac{h_+^2 - h_+h_- + h_-^2}{12}\;u^{(4)}\\
  \ -\ \frac{(h_+ - h_-)(h_+^2 + h_-^2)}{60}\;u^{(5)}\ +\ \O(h^4)\,,
\end{multline}
where we used the fact that $\O\,(h_-^4)\ =\ \O\,(h_+^4)\ =\ \O\,(h^4)\,$, which follows from \eqref{eq:jac}. Since $u\,(x_{\,j})$ is the solution of \eqref{eq:bvp}, then $-u''\ +\ \lambda^2\,u\ \equiv\ 0\,$.

By assuming that the mapping $x(q)$ is sufficiently smooth, we obtain
\begin{align*}
  h_+\ &=\ x\,(q_{\,j+1})\ -\ x\,(q_{\,j})\ =\ h\,x_{\,q}\ +\ \frac{h^2}{2}\; x_{\,qq}\ +\ \frac{h^3}{6}\;x_{\,qqq}\ +\ \O\,(h^4)\,, \\
  h_-\ &=\ x\,(q_{\,j})\ -\ x\,(q_{\,j-1})\ =\ h\,x_{\,q}\ -\ \frac{h^2}{2}\; x_{\,qq}\ +\ \frac{h^3}{6}\;x_{\,qqq}\ +\ \O\,(h^4)\,.
\end{align*}
The last two expansions allow us to estimate asymptotically the following combinations, which appear as coefficients in \eqref{eq:psij}:
\begin{align*}
  h_+\ -\ h_-\ &=\ h^2\,x_{qq}\ +\ \O\,(h^4)\,, \\
  h_+^2\ -\ h_+h_-\ +\ h_-^2\ &=\ h^2\,x_q^2\ +\ \O\,(h^4)\,, \\
  (h_+\ -\ h_-)\,(h_+^2\ +\ h_-^2)\ &=\ \O\,(h^4)\,.
\end{align*}
By substituting these estimations into \eqref{eq:psij}, one can conclude
\begin{equation*}
  \psi_{\,j}\ =\ -\frac{h^2}{3}\;\Bigl[x_{\,qq}\,u_{\,xxx}\ +\ \frac{1}{4}\;x_{\,q}^2\,u_{\,xxxx}\Bigr]\ +\ \O(h^4)\,.
\end{equation*}
Thus, it shows that in general the scheme \eqref{eq:nonuni} is second order accurate. However, one can notice that the scheme will be exceptionally the fourth order accurate if the mapping $x\,(q)$ satisfies the following equation:
\begin{equation}\label{eq:fourth}
  x_{\,qq}\, u_{\,xxx}\ +\ \frac{1}{4}\;x_{\,q}^{\,2}\,u_{\,xxxx}\ \equiv\ 
  x_{\,qq}\, u_{\,xxx}\ +\ \frac{1}{4}\;x_{\,q}\,u_{\,xxxq}\ =\ 0\,,
\end{equation}
where $u\,(x)$ is a solution to equation \eqref{eq:bvp}. Finally, we can rewrite equation \eqref{eq:fourth} in a compact form
\begin{equation}\label{eq:last}
  (u_{\,xxx})^{\,3/4}\,\bigl[\,(u_{\,xxx})^{\,1/4}\,x_{\,q}\bigr]_{\,q}\ =\ 0\,.
\end{equation}
After noticing that $u_{\,x\,x\,x} \propto u_{\,x}$ thanks to \eqref{eq:bvp}, equation \eqref{eq:last} can be recast into an equivalent, but more familiar form:
\begin{equation*}
  \bigl[\,(u_{\,x})^{1/4}\,x_{\,q}\bigr]_{\,q}\ =\ 0\,,
\end{equation*}
which is to be compared with \eqref{eq:elliptic} and \eqref{eq:example}. This observation demystifies the surprising performance of $\beta\ =\ \frac{1}{4}\,$. It seems that this observation was made for the first time by \textsc{Degtyarev} \etal \cite{Degtyarev1987a} in a slightly different form. In order to make this argument rigorous, one needs also to prove the stability of the scheme. It is done below in Section~\ref{sec:stab}.

\begin{remark}\label{rem:1}
\textsc{Degtyarev} \etal \cite{Degtyarev1987a} also constructed the `optimal' monitoring function by performing the local error $\psi_j$ analysis. As a result, they obtain a similar monitoring function:
\begin{equation}\label{eq:degt}
  \omega\,(x)\ =\ \bigl(\eps\,\abs{\,u_{\,x\,x\,x}\,}\,\bigr)^{\frac{1}{4}}\,.
\end{equation}
The same monitoring function was used in \cite{Ferreira1993} as well.

First of all, we note that the third derivative in this problem is positive and the absolute value can be omitted. Moreover, thanks to Lemma~\ref{lem:prop} the resulting meshes coincide. However, there is a clear advantage of the monitoring function $\omega\ =\ \bigl(\,u_{\,x}\,\bigr)^{\,\frac{1}{4}}$ proposed in our study over \eqref{eq:degt} at the level of the numerical implementation, which were not performed in \cite{Degtyarev1987a, Ferreira1993}. Namely, the approximation of the third order derivatives on non-uniform meshes is not a trivial task, especially near the boundaries. The authors of \cite{Ferreira1993} were aware of this problem and they suggested to use finite difference formulas developed in \cite{Garcia-Archilla1991}. These difficulties can be completely avoided by lowering the order of derivatives in the monitoring function by using the governing equation \eqref{eq:bvp}.

As a result, we can issue another practical recommendation: the order of derivatives in the monitoring function should be lowered as much as possible by using all available information about the solutions.
\end{remark}


\subsection{Remarks on the practical implementation}

Above we performed an analysis of the finite difference approximation error which gave us an `optimal' monitoring function\footnote{Taking into account the Remark~\ref{rem:1} we know that the monitoring function \eqref{eq:degt} produces the same grid as the function $\omega\,(x)\ =\ \bigl(\,u_{\,x}\,\bigr)^{\frac{1}{4}}\,$.} $\omega\,(x)\ =\ \bigl(\,u_{\,x\,x\,x}\,\,\bigr)^{\frac{1}{4}}\ \propto\ \bigl(\,u_{\,x}\,\bigr)^{\frac{1}{4}}$ yielding to the fourth order supraconvergence phenomenon. However, there might be some practical implementation problems with this choice \eqref{eq:degt} of the monitoring function. Indeed, on solutions to BVP \eqref{eq:bvp} this function takes almost zero values on large portions of the computational $(0,\,\ell)$ starting from the point $x\ =\ 0\,$. In practice it can result in a very few nodes placed on this sub-domain, which is not desirable. To correct this shortcoming, we can use instead another monitoring function, which is bounded from below by a positive constant. This is why in our numerical computations we employed also the monitoring function \eqref{eq:moni}. In this Section we consider the following `optimal' sub-family of functions \eqref{eq:moni}, which are bounded from below by $1\,$:
\begin{equation}\label{eq:alpha}
  \omega_{\,\alpha}\,(x)\ =\ 1\ +\ \alpha\;\abs{\,u_{\,x}}^{\,\frac{1}{4}}\,.
\end{equation}
The question one might ask is how different are the grids produced by monitoring functions $\omega\,(x)\ =\ \bigl(\,u_{\,x}\,\bigr)^{\frac{1}{4}}$ and $\omega_{\,\alpha}\,(x)$ defined in \eqref{eq:alpha}? It turns out that the answer on this question does not depend on the exponent $\beta\ =\ \frac{1}{4}$ and we show a more general
\begin{lemma}\label{lem:2}
Two grids $x\ =\ x_{\,\beta}\,(q)$ and $x\ =\ x_{\,\beta}^{\,\alpha}\,(q)$ produced by monitoring functions $\omega_{\,\beta}\,(x)\ =\ \bigl(\,u_{\,x}\,\bigr)^{\,\beta}$ and \eqref{eq:moni} correspondingly are close enough provided that $\alpha\ >\ 0$ is large, \ie
\begin{equation}\label{eq:est}
  \lim_{\alpha\ \to\ +\infty} \norm{x_{\,\beta}\ -\ x_{\,\beta}^{\,\alpha}}_{\infty}\ =\ 0\,.
\end{equation}
\end{lemma}

\begin{proof}
By integrating once the elliptic equation \eqref{eq:elliptic} it follows that
\begin{equation}\label{eq:one}
  \omega_{\,\beta}\,(x)\,x_{\,\beta,\,q}\ \equiv\ \frac{\lambda^{\,\beta\ -\ 1}}{\beta}\;\bigl(1\ -\ \ue^{-\lambda\,\beta\,\ell}\bigr) \ \defeq\ C\,.
\end{equation}
For the second mapping we integrate once the equation
\begin{equation*}
  \od{}{q}\;\biggl[\omega_{\,\alpha,\,\beta}\,(x)\,\od{x_{\,\beta}^{\,\alpha}}{q}\biggr]\ =\ 0\,,
\end{equation*}
yielding to
\begin{equation*}
  \omega_{\,\alpha,\,\beta}\,(x)\,x_{\,\beta,\,q}^{\,\alpha}\ =\ \ell\ +\ \alpha\,C\,,
\end{equation*}
or equivalently
\begin{equation}\label{eq:two}
  \bigl(1\ +\ \alpha\,\omega_{\,\beta}\,(x)\bigr)\,x_{\,\beta,\,q}^{\,\alpha}\ =\ \ell\ +\ \alpha\,C\,.
\end{equation}
From equations \eqref{eq:one} and \eqref{eq:two} we can find:
\begin{equation*}
  x_{\,\beta,\,q}\ =\ \frac{C}{\omega_{\,\beta}\,(x)}\,, \qquad
  x_{\,\beta,\,q}^{\,\alpha}\ =\ \frac{\ell\ +\ \alpha\,C}{1\ +\ \alpha\,\omega_{\,\beta}\,(x)}\,.
\end{equation*}
By taking the difference between two last equations we obtain:
\begin{equation*}
  x_{\,\beta,\,q}^{\,\alpha}\ -\ x_{\,\beta,\,q}\ =\ \frac{\ell\,\omega_{\,\beta}\,(x)\ -\ C}{\omega_{\,\beta}\,(x)\,\bigl(1\ +\ \alpha\,\omega_{\,\beta}\,(x)\bigr)}\,.
\end{equation*}
In the right hand side the monitoring function $\omega_{\,\beta}\,(x)$ is bounded and together with the constant $C$ does not depend on $\alpha\,$. Thus, for $\alpha\ \to\ +\infty$ we have a uniform estimate:
\begin{equation*}
  \norm{x_{\,\beta,\,q}\ -\ x_{\,\beta,\,q}^{\,\alpha}}_{\infty}\ =\ \O\,\Bigl(\frac{1}{\alpha}\Bigr)\,.
\end{equation*}
Since in boundary nodes $q\ =\ 0$ and $q\ =\ 1$ the difference $x_{\,\beta}\ -\ x_{\,\beta}^{\,\alpha}$ vanishes. Thus, from the last estimate on the derivatives difference and by applying the embedding theorem\footnote{Let us remind the embedding theorems in the simplest 1D case needed for our purposes. Let us take a continuously differentiable function $u\,(x)\ \in\ C^{\,1}\,[\,0,\,\ell\,]$, which takes zero values on the boundaries $x\ =\ 0$ and $x\ =\ \ell\,$. Then, we have two inequalities:
\begin{equation*}
  \norm{\,u\,}_{\,\infty}\ \leq\ \frac{\sqrt{\ell}}{2}\;\norm{\,u_{\,x}\,}_{\,L_2}\,, \qquad \norm{\,u_{\,x}\,}_{\,L_2}\ \leq\ \sqrt{\ell}\,\norm{\,u_{\,x}\,}_{\,\infty}\,,
\end{equation*}
which yield the required estimation \eqref{eq:esti}.} it follows that
\begin{equation}\label{eq:esti}
  \norm{x_{\,\beta}\ -\ x_{\,\beta}^{\,\alpha}}_{\infty}\ =\ \O\,\Bigl(\frac{1}{\alpha}\Bigr)\,.
\end{equation}
By taking the limit $\alpha\ \to\ +\infty$ we obtain the required statement \eqref{eq:est}.
\end{proof}

\begin{remark}
We would like to mention that for $\beta\ =\ \frac{1}{4}$ the grid $x\ =\ x_{\,\beta}^{\,\alpha}\,(q)$ will not be optimal anymore as it was illustrated in our numerical simulations. However, the real question which should be asked for a general value of the parameter $\beta\,$: how accurate are solutions obtained on close grids $x\ =\ x_{\,\beta}\,(q)$ and $x\ =\ x_{\,\beta}^{\,\alpha}\,(q)\,$. We do not have a theoretical answer on this question. However, the numerical results reported in Table~\ref{tab:res} show that for a fixed $\beta$ and large values of $\alpha\ \simeq\ 10^{\,3}\,\ldots\,10^{\,4}$ the error depends only marginally on $\alpha\,$, which is consistent with Lemma~\ref{lem:2}.
\end{remark}


\section{Stability and convergence}
\label{sec:stab}

Above we showed that the approximation error of the equation on the `optimal' grid is of order four in $h\,$. However, it is more important to know whether the convergence is of order four? In this Section we show that it is indeed the case.

In order to show this result, we first recall a technical Lemma from \cite{Samarskii2001}:
\begin{lemma}\label{lem:sam}
Consider the following difference scheme:
\begin{equation}\label{eq:star}
  a_{\,j}\,u_{\,j-1}\ +\ c_{\,j}\,u_{\,j}\ +\ b_{\,j}\,u_{\,j+1}\ =\ f_{\,j}\,, \qquad
  j\ =\ 1,\,2,\,\ldots,\,N\,-\,1\,,
\end{equation}
together with boundary conditions:
\begin{equation*}
  u_{\,0}\ =\ u\,(0)\,, \qquad
  u_{\,N}\ =\ u\,(\ell)\,.
\end{equation*}
If coefficients $\bigl\{a_{\,j},\,b_{\,j},\,c_{\,j}\bigr\}_{\,j\,=\,1}^{\,N\,-\,1}$ satisfy the condition
\begin{equation}\label{eq:conda}
  \abs{\,c_{\,j}\,}\ \geq\ \abs{\,a_{\,j}\,}\ +\ \abs{\,b_{\,j}\,}\ +\ \delta\,, \qquad \delta\ >\ 0\,,
\end{equation}
then, the difference problem \eqref{eq:star} possesses the unique solution for any data $u\,(0)\,$, $u\,(\ell)$ and $\bigl\{f_{\,j}\bigr\}_{\,j\,=\,1}^{\,N\,-\,1}$ and we have an estimate\footnote{We remind that in a finite dimensional space (\ie when $h$ is fixed) all norms are equivalent. Thus, the same estimate holds in other norms modulo a constant factor.}:
\begin{equation}\label{eq:lest}
  \norm{u_{\,h}}_{\,\infty}\ \leq\ \max\Bigl\{\;\abs{\,u\,(0)\,},\;\abs{\,u\,(\ell)\,},\;\frac{1}{\delta}\cdot\max_{1\,\leq\,j\,\leq\,N\,-\,1}\abs{\,f_{\,j}\,}\;\Bigr\}\,,
\end{equation}
where $u_{\,h}$ denotes the discrete grid function.
\end{lemma}

\begin{proof}
See \cite{Samarskii2001}.
\end{proof}

This Lemma will be used during the proof of the main
\begin{theorem}
The finite difference solution to problem \eqref{eq:nonuni} constructed using the `optimal' grid is stable and converges to the exact solution \eqref{eq:sol} with order $\O\,(h^{\,4})\,$.
\end{theorem}

\begin{proof}
Our problem \eqref{eq:bvp} and, thus, its finite difference discretization \eqref{eq:nonuni} are homogeneous (\ie the right hand side is zero). However, we can prove our result in a more general case, where we have an `exterior forcing' $\bigl\{f_{\,j}\bigr\}_{\,j\,=\,1}^{\,N\,-\,1}$ and arbitrary boundary conditions. Henceforth, instead of \eqref{eq:nonuni} we consider a more general difference problem:
\begin{equation}\label{eq:prob}
  -\frac{1}{h_j}\;\biggl[\,\frac{u_{j+1}\ -\ u_j}{h_{j+1/2}}\ -\ \frac{u_j\ -\ u_{j-1}}{h_{j-1/2}}\,\biggr]\ +\ \lambda^2\,u_j\ =\ f_{\,j}\,, 
  \qquad j\ =\ 1,\,2,\, \ldots,\, N\, -\, 1\,,
\end{equation}
together with boundary conditions:
\begin{equation}\label{eq:bc}
  u_{\,0}\ =\ \mu_{\,0}\,, \qquad
  u_{\,N}\ =\ \mu_{\,\ell}\,.
\end{equation}
First we establish the stability property of this scheme in $l_{\,\infty}$ norm:
\begin{equation*}
  \norm{\,u_{\,h}\,}_{\,\infty}\ \eqdef\ \max_{0\,\leq\,j\,\leq\,N}\,\abs{\,u_{\,j}\,}\,.
\end{equation*}
The problem \eqref{eq:prob}, \eqref{eq:bc} can be written in the operator short-hand form:
\begin{equation*}
  \Ll^h\, u_h\ =\ f_{\,h}\,.
\end{equation*}
The norm of the right hand side will be defined as
\begin{equation*}
  \norm{f_{\,h}}_{\,\infty}^{\,\prime}\ \eqdef\ \max\Bigl\{\,\abs{\,\mu_{\,0}\,},\;\abs{\,\mu_{\,\ell}\,},\;\max_{1\,\leq\,j\,\leq\,N\,-\,1}\,\abs{\,f_{\,j}\,}\,\Bigr\}\,.
\end{equation*}
We seek to apply Lemma~\ref{lem:sam} to finite difference problem \eqref{eq:prob}, \eqref{eq:bc}. Indeed, the coefficients $\bigl\{a_{\,j},\,b_{\,j},\,c_{\,j}\bigr\}_{\,j\,=\,1}^{\,N\,-\,1}$ can be easily computed:
\begin{align*}
  a_{\,j}\ &=\ -\,\frac{1}{h_{\,j}\,h_{\,j-\frac{1}{2}}}\ <\ 0\,, \\
  b_{\,j}\ &=\ -\,\frac{1}{h_{\,j}\,h_{\,j+\frac{1}{2}}}\ <\ 0\,, \\
  c_{\,j}\ &=\ -\,a_{\,j}\ -\ b_{\,j}\ +\ \lambda^{\,2}\ >\ 0\,.
\end{align*}
Thus, condition \eqref{eq:conda} is satisfied with $\delta\ \eqdef\ \lambda^{\,2}\ >\ 0\,$. Moreover, Lemma~\ref{lem:sam} gives us the following estimate of the solution:
\begin{equation*}
  \norm{\,u_{\,h}\,}_{\,\infty}\ \leq\ \max\Bigl\{\;\abs{\,\mu_{\,0}\,},\;\abs{\,\mu_{\,\ell}\,},\;\frac{1}{\lambda^{\,2}}\cdot\max_{1\,\leq\,j\,\leq\,N\,-\,1}\abs{\,f_{\,j}\,}\;\Bigr\}\,.
\end{equation*}
In other words we have
\begin{equation*}
  \norm{\,u_{\,h}\,}_{\,\infty}\ \leq\ C_{\,s}\,\norm{f_{\,h}}_{\,\infty}^{\,\prime}\,,
\end{equation*}
where we introduced a constant $C_{\,s}\ \eqdef\ \max\Bigl\{\,1,\,\dfrac{1}{\lambda^{\,2}}\,\Bigr\}\,$. Thus, we just showed the numerical solution stability in $l_{\,\infty}$ norm. Earlier in Section~\ref{sec:expla}, we established the approximation property and its order ($\O\,(h^{\,4})$ on the `optimal' grid). Now, by applying the \textsc{Lax}--\textsc{Richtmyer} fundamental theorem of (linear) numerical analysis \cite{Lax1956}, we conclude the proof of our Theorem.
\end{proof}

\begin{remark}
Notice that the `optimal' grid is constructed in practice by solving numerically \eqref{eq:nbvp} even if exact solutions are available in some particular cases. In order to enjoy the fourth order accuracy, one has to solve numerically the elliptic equation \eqref{eq:elliptic} at least to the second order accuracy.
\end{remark}


\section{Discussion}
\label{sec:concl}

In the present note we provided an example of a boundary value problem, which exhibits the boundary layer phenomenon, as a prototype of more complex problems arising \eg in Fluid Mechanics \cite{Schlichting2000}. When this equation is discretized with central finite differences on a uniform grid, one formally obtains the second order accuracy uniformly in space. However, the proportionality constant in the consistency error term becomes unacceptably large in the boundary layer. So, this solution cannot be satisfactory for more complex 3D problems. Consequently, we proposed to keep precisely the same scheme, but to modify the grid using the equidistribution principle. As a result we show first that the proportionality constant becomes quasi-uniform in space (no substantial increase in the boundary layer), but the main gain is in the order of convergence, which becomes equal to \emph{four} provided that the mesh nodes are distributed accordingly. One can even ask a question in view of the study \cite{Chhay2011}: is it better to redistribute the nodes in order to improve substantially the scheme accuracy instead of preserving some symmetries at the discrete level? The numerical results of the present study seem to favour the former possibility.

This example is very simple and instructive since it shows how poor is our current understanding of the numerical analysis on non-uniform adaptive\footnote{The grid becomes \emph{moving} for time-dependent problems. See \cite{Khakimzyanov2015a} for some examples on hyperbolic problems.} grids. Similar techniques can be used to solve more realistic problems involving boundary layers, blow up phenomena \cite{Budd2001} and other rapid variations in space/time of the solution. Recently we employed these techniques to solve numerically hyperbolic conservation laws \cite{Khakimzyanov2015a}. From our experience it follows that one can expect a substantial improvement of the scheme properties when adaptive grids are \emph{carefully} employed. However, the negative examples we provide show also that the application of adaptive grids does not necessarily lead \emph{per se} to the increase in the accuracy. Consequently, the critical analysis of obtained numerical results is always needed.

The perspectives opened by this study include the description of \emph{optimal} monitoring functions for some classes of problems. In the absence of such optimal choice, one has to describe at least the values of free parameters (\eg $\alpha$, $\beta$) which yield converged numerical solutions with desired properties. Here we made a proposition for a singularly perturbed Sturm--Liouville problem. However, the variety of boundary layers encountered in practice is much richer. It looks like a successful numerics can be achieved only in conjunction with deep analytical understanding of the problem in hands.


\subsection*{Acknowledgments}
\addcontentsline{toc}{subsection}{Acknowledgments}

This research was supported by RSCF project No 14-17-00219. D.~\textsc{Dutykh} acknowledges the support of the CNRS under the PEPS InPhyNiTi project FARA and project \No EDC26179 --- ``\textit{Wave interaction with an obstacle}'' as well as the hospitality of the Institute of Computational Technologies SB RAS during his visit in October 2015. The authors would also like to thank Prof. Laurent \textsc{Gosse} (CNR, Italy) for stimulating discussions on singularly perturbed Sturm--Liouville problems.


\appendix
\section{A first order problem}
\label{sec:app}

In this Appendix we apply the technique developed above to a first order problem to show the handling of first order derivatives and Initial Value Problems (IVP). Consider the following \textsc{Cauchy}-type problem:
\begin{equation}\label{eq:first}
  \od{u}{x}\ +\ \lambda\,u\ =\ f\,(x)\,, \qquad 0\ <\ x\ \leq\ \ell\,, \qquad u\,(0)\ =\ \mu_{\,0}\,,
\end{equation}
where $\lambda\ =\ \const\ >\ 0\,$. Notice that a similar (but simpler) problem was considered in \cite[Section~\textsection2]{Ferreira1993}. However, the Authors of \cite{Ferreira1993} ill posed the difference problem, which resulted in its non-solvability for every even number of discretization points $N\ =\ 2\,N_{\,0}\,$, $N_{\,0}\ \in\ \N\,$. This drawback will be corrected here.

We approximate the continuous problem \eqref{eq:first} by the following difference problem (on a uniform grid for the moment):
\begin{equation}\label{eq:fduni}
  \frac{u_{\,j+1}\ -\ u_{\,j-1}}{2\,\Delta x}\ +\ \lambda\,u_{\,j}\ =\ f_{\,j}\,, \qquad
  j\ =\ 1,\,2,\,\ldots,\,N\,-\,1\,,
\end{equation}
together with initial conditions:
\begin{equation}\label{eq:ic2}
  u_{\,0}\ =\ \mu_{\,0}\,, \qquad u_{\,1}\ =\ \mu_{\,1}\,.
\end{equation}
The crucial point here is to pose \emph{two} initial conditions, while the Authors of \cite{Ferreira1993} pose only one by analogy with the continuous problem. For difference problems the number of initial conditions necessary to have the well posedness usually depends on the degree of the characteristic equation induced by this scheme. For instance, for difference problem \eqref{eq:fduni}, the characteristic equation is of the second degree. Thus, two initial conditions are needed in nodes $x\ =\ x_{\,0}\ \equiv\ 0$ and $x\ =\ x_{\,1}\ \equiv\ \Delta x\,$.

One might ask the question is how to compute the value $\mu_{\,1}\,$, which is not provided by the continuous formulation? There are various approaches. For instance, one might integrate IVP \eqref{eq:first} using a \textsc{Runge}--\textsc{Kutta} (of at least second order accuracy) method \cite{Press2007}, for example. In our study we suggest to use a local \textsc{Taylor} expansion together with the information provided by the governing equation:
\begin{multline*}
  \mu_{\,1}\ \equiv\ u\,(x_{\,1})\ \stackrel{\displaystyle{\mathrm{\textsc{Taylor}}}}{=}\ \underbrace{u\,(x_{\,0})}_{\displaystyle{\equiv\ \mu_{\,0}}}\ +\ u^{\,\prime}\,(x_{\,0})\cdot\Delta x\ +\ \O\,(\Delta x^{\,2})\ =\\ 
  \mu_{\,0}\ +\ \underbrace{\bigl(f\,(x_{\,0})\ -\ \lambda\,u\,(x_{\,0})\bigr)}_{\displaystyle{\eqref{eq:first}}}\,\Delta x\ +\ \O\,(\Delta x^{\,2})\ =\\
  \mu_{\,0}\,\bigl(1\ -\ \lambda\,\Delta x\bigr)\ +\ f_{\,0}\,\Delta x\ +\ \O\,(\Delta x^{\,2})\,.
\end{multline*}
Thus, after neglecting higher order terms, we recover the missing data:
\begin{equation*}
  \mu_{\,1}\ \eqdef\ \mu_{\,0}\,\bigl(1\ -\ \lambda\,\Delta x\bigr)\ +\ f_{\,0}\,\Delta x\,.
\end{equation*}
Moreover, in this way we conserve the second approximation order uniformly.

Now, we consider problem \eqref{eq:first} posed on the reference domain $Q\,$:
\begin{equation}\label{eq:23}
  \frac{1}{J}\;\od{v}{q}\ +\ \lambda\,v\ =\ f\,(q)\,, \qquad 0\ <\ q\ \leq\ 1\,,
  \qquad v\,(0)\ =\ \mu_{\,0}\,,
\end{equation}
where $J$ is the \textsc{Jacobian} of the underlying coordinate transformation $x\ =\ x\,(q)\,$. The last problem can be discretized on a uniform grid $Q_{\,h}$ with step $h$ which covers $Q\,$. The last scheme \eqref{eq:23} can be rewritten as a finite difference scheme on a non-uniform grid in the physical domain $\I_{\,h}\,$:
\begin{equation}\label{eq:24}
  \frac{u_{\,j+1}\ -\ u_{\,j-1}}{h_{\,j+1/2}\ +\ h_{\,j-1/2}}\ +\ \lambda\,u_{\,j}\ =\ f_{\,j}\,, \qquad
  j\ =\ 1,\,2,\,\ldots,\,N\,-\,1\,,
\end{equation}
with initial conditions \eqref{eq:ic2}. We can study the approximation properties of the scheme \eqref{eq:24} using similar methods to what we did earlier in Section~\ref{sec:expla}:
\begin{equation*}
  \psi_{\,j}\ =\ \frac{h^{\,2}}{2}\;\bigl(u_{\,x\,x}\bigr)^{\,\frac{2}{3}}\,\Bigl[\,x_{\,q}\,\bigl(u_{\,x\,x}\bigr)^{\,\frac{1}{3}}\,\Bigr]_{\,q}\ +\ \O\,(h^{\,4})\,.
\end{equation*}
From the last equation it follows even on the non-uniform grid we still have the second order approximation to the differential equation \eqref{eq:first}. However, there exists a special transformation $x\ =\ x\,(q)\,$, which can be computed with the equidistribution method for a special choice of the monitoring function $\omega\,(x)\ =\ \bigl(u_{\,x\,x}\bigr)^{\,\frac{1}{3}}\,$. On this grid one enjoys the supraconvergence phenomenon with the fourth approximation order. It can be shown using different methods that on this grid the solution is stable and, thus, it converges to the true solution with the fourth order as well (by application of the \textsc{Lax}--\textsc{Richtmyer} equivalence theorem). We note that the monitoring function $\omega\,(x)\ =\ \bigl(u_{\,x\,x}\bigr)^{\,\frac{1}{3}}$ was obtained in \cite{Ferreira1993} as well.

We would like to mention that the `optimal' monitoring function contains the derivatives higher than the order of the problem. Thus, according to Remark~\ref{rem:1} we should better lower this order for practical implementations. It can be done by differentiating once the governing equation \eqref{eq:first}:
\begin{equation*}
  u_{\,x\,x}\ \equiv\ f_{\,x}\ -\ \lambda\,f\ +\ \lambda^{\,2}\,u\,.
\end{equation*}
Thus, the same `optimal' grid can be obtained using another monitoring function, which contains lower order derivatives. Numerically this observation is quite useful.

As a final remark we have to say that if one uses the `optimal' grid, the missing initial data $\mu_{\,1}$ has to be computed to the fourth order accuracy as well. To achieve this goal one may use again the local \textsc{Taylor} expansion (to the 4\up{th} order this time) or perform one step of the 4\up{th} order \textsc{Runge}--\textsc{Kutta} scheme \cite{Press2007}.


\subsection{A simple analytical case}

For testing purposes of the numerical algorithms we consider here the homogeneous case $f\,(x)\ \equiv\ 0\,$, where everything can be constructed analytically, including the exact solution to measure the error of the numerical one. So, the exact solution to \textsc{Cauchy}-type problem \eqref{eq:first} is
\begin{equation*}
  u\,(x)\ =\ \mu_{\,0}\,\ue^{-\,\lambda\,x}\,.
\end{equation*}
The exact monitoring function is proportional to
\begin{equation*}
  \omega\,(x)\ \propto\ \ue^{-\,\lambda\,x/3}\,.
\end{equation*}
It is not difficult to check that for this monitoring function, the equidistribution method generates the following transformation:
\begin{equation*}
  x\,(q)\ =\ -\,\frac{3}{\lambda}\;\ln\bigl[\,1\ -\ q\ +\ q\,\ue^{-\,\lambda\,\ell/3}\,\bigr]\,,
\end{equation*}
which produces the `optimal' non-uniform grid. Finally, we compute also exactly the missing initial data $\mu_{\,1}\,$. Our method still relies on the 4\up{th} order local \textsc{Taylor} expansion:
\begin{equation*}
  u\,(x_{\,1})\ =\ u\,(x_{\,0})\ +\ u^{\,\prime}\,(x_{\,0})\cdot h_{\,1/2}\ +\ u^{\,\dprime}\,(x_{\,0})\;\frac{h_{\,1/2}^{\,2}}{2}\ +\ u^{\,\trprime}\,(x_{\,0})\;\frac{h_{\,1/2}^{\,3}}{6}\ +\ \O\,(h_{\,1/2}^{\,4})\,.
\end{equation*}
From equation \eqref{eq:first}, taking into account that $f\,(x)\ \equiv\ 0\,$, we have
\begin{equation*}
  u^{\,\prime}\,(x_{\,0})\ =\ -\lambda\,\mu_{\,0}\,, \qquad
  u^{\,\dprime}\,(x_{\,0})\ =\ \lambda^{\,2}\,\mu_{\,0}\,, \qquad
  u^{\,\trprime}\,(x_{\,0})\ =\ -\lambda^{\,3}\,\mu_{\,0}\,.
\end{equation*}
Henceforth, we can reconstruct the missing initial data to the fourth order accuracy as follows:
\begin{equation*}
  \mu_{\,1}\ =\ \mu_{\,0}\cdot\biggl[\,1\ -\ \lambda\,h_{\,1/2}\ +\ \frac{\lambda^{\,2}\,h_{\,1/2}^{\,2}}{2}\ -\ \frac{\lambda^{\,3}\,h_{\,1/2}^{\,3}}{6}\,\biggr]\,.
\end{equation*}
This concludes the presentation of another example of successful application of the equidistribution method to increase significantly the order of the underlying scheme.


\bigskip
\addcontentsline{toc}{section}{References}
\bibliographystyle{abbrv}

\bigskip

\end{document}